\documentclass{amsart}

\usepackage{amsmath,amssymb,amsthm}

\newtheorem{prop}{Proposition}
\newtheorem{thm}[prop]{Theorem}
\newtheorem{lem}[prop]{Lemma}

\theoremstyle{definition}

\def\co{\colon\thinspace}

\newcommand{\CP}{\mathbb{C}\mathrm{P}}

\newcommand{\rmd}{\mathrm{d}}

\newcommand{\tlambda}{\tilde{\lambda}}

\newcommand{\omegacan}{\omega_{\mathrm{can}}}
\newcommand{\omegamag}{\omega_{\mathrm{mag}}}

\newcommand{\bfp}{\mathbf{p}}

\newcommand{\bfq}{\mathbf{q}}

\newcommand{\R}{\mathbb{R}}
\newcommand{\RP}{\mathbb{R}\mathrm{P}}

\newcommand{\bft}{\mathbf{t}}

\newcommand{\bfu}{\mathbf{u}}

\newcommand{\xican}{\xi_{\mathrm{can}}}

\newcommand{\Z}{\mathbb Z}

\begin{document}

\author[H. Geiges]{Hansj\"org Geiges}
\address{Mathematisches Institut, Universit\"at zu K\"oln,
Weyertal 86--90, 50931 K\"oln, Germany}
\email{geiges@math.uni-koeln.de}
\author[K. Zehmisch]{Kai Zehmisch}
\address{Fakult\"at f\"ur Mathematik, Ruhr-Universit\"at Bochum,
Universit\"atsstra{\ss}e 150, 44780 Bochum, Germany}
\email{kai.zehmisch@rub.de}

\title[Symplectic fillings of unit cotangent bundles]{Symplectic fillings of
unit cotangent bundles of hyperbolic surfaces}

\date{}

\begin{abstract}
We consider strong symplectic fillings
of the unit cotangent bundle of a hyperbolic surface,
equipped with its canonical contact structure.
We show that every finitely presentable group can be realised as
the fundamental group of such a filling.
\end{abstract}

\thanks{This research is part of a project in the SFB/TRR 191
\textit{Symplectic Structures in Geometry, Algebra and Dynamics}, 
funded by the DFG (Project-ID 281071066 -- TRR 191)}

\maketitle

%%%%%%%%%%%%%%%%%%%%%%%%%%%%%%%%%%%%%%%%

\section{Introduction}
Much of the literature on symplectic fillings of contact manifolds
has been concerned with establishing uniqueness results of various
kinds, be it up to symplectomorphism, deformation equivalence,
or diffeomorphism of the filling. The exemplar of such results
is the unique strong symplectic fillability of the
sphere $S^{2n-1}$ with its standard contact structure by
the standard symplectic ball --- up to symplectomorphism for
$n=2$ (Gromov~\cite[p.~311]{grom85}
and McDuff~\cite[Theorem~1.7]{mcdu90}); up to diffeomorphism
for $n\geq 3$ (Eliashberg--Floer--McDuff; see \cite[Theorem~1.5]{mcdu91}).
Here the filling is understood to be symplectically aspherical
for $n\geq 3$, that is, the cohomology class of the symplectic form
evaluates to zero on homology classes represented by spheres;
for $n=2$ it suffices to demand minimality of the
filling, that is, the absence of symplectically
embedded $2$-spheres of self-intersection $-1$. This requirement
prevents the obvious modification of the filling by symplectic blow-ups.

For some recent uniqueness results we refer the reader to
\cite{bgz19} and \cite{gkz23}; the introductions to
these papers contain surveys of the literature
on the subject.

In the present paper, by contrast, we describe
ways to modify the diffeomorphism or symplectomorphism type
of a given filling. Specifically, we apply these modifications to
strong symplectic fillings of
unit cotangent bundles $ST^*\Sigma$ of closed, oriented surfaces
$\Sigma=\Sigma_g$ of genus $g\geq2$, equipped with
the canonical contact structure~$\xican$. We shall be referring to these
surfaces as \emph{hyperbolic surfaces}, even though the choice of
Riemannian metric is irrelevant.

For lower genus, or for more restrictive types of filling
(cf.~\cite[Chapter~5]{geig08}), there are once again a number of
uniqueness statements. In the following list
it is understood that all unit cotangent
bundles are equipped with the contact structure~$\xican$:
\begin{itemize}
\item[-] For $ST^*S^2$, the unit disc bundle (with its canonical symplectic
form~$\omegacan$) is the unique minimal strong symplectic filling
up to diffeomorphism, and up to symplectomorphism if one fixes the cohomology
class of the symplectic form (McDuff~\cite{mcdu90}).
\item[-] Stein fillings of $ST^*S^2$ are unique up to Stein
homotopy (Hind~\cite{hind00}).
\item[-] All Stein fillings of $ST^*T^2$ are homeomorphic
(Stipsicz~\cite{stip02}).
\item[-] All minimal strong symplectic fillings of $ST^*T^2$ are
deformation equivalent, and in particular diffeomorphic
(Wendl~\cite{wend10}).
\item[-] For hyperbolic surfaces~$\Sigma$, all Stein fillings
are $s$-cobordant rel boundary, and exact symplectic fillings satisfy
homological restrictions (Sivek and Van Horn-Morris~\cite{siva17},
Li--Mak--Yasui~\cite{lmy17}).
\end{itemize}

For hyperbolic surfaces $\Sigma$, however, non-diffeomorphic
strong symplectic fillings of $ST^*\Sigma$ \emph{do} exist, as was
first observed by Li--Mak--Yasui~\cite[Proof of Corollary~1.6]{lmy17}.
One starts with the $4$-manifold $[-1,1]\times ST^*\Sigma$, equipped
with a symplectic form such that both boundary components are strongly
convex; in particular, they carry an induced contact structure compatible
with the
boundary orientation, and the contact structure on $\{1\}\times
ST^*\Sigma$ equals~$\xican$. Such a symplectic form has been
constructed in \cite{mcdu91} and \cite{geig95}. One then caps
off the boundary component $\{-1\}\times ST^*\Sigma$ by the method of
Eliashberg~\cite{elia04} and Etnyre~\cite{etny04} to obtain a minimal strong
symplectic filling
of $(ST^*\Sigma,\xican)$ with Betti number $b_2^+$ arbitrarily
large. In fact, since we need only cap off a \emph{strongly}
convex boundary, the earlier result of Etnyre--Honda~\cite{etho02}
can be used to construct such a cap; see also \cite{geig06}.

In the present note we construct infinitely many non-diffeomorphic
symplectic fillings of $(ST^*\Sigma,\xican)$ by arguably more
elementary means. Our construction has the added advantage of
giving us control over the fundamental group.

\begin{thm}
\label{thm:main}
Every finitely presentable group can be realised as the fundamental
group of a strong symplectic filling of $(ST^*\Sigma,\xican)$.
\end{thm}

Starting point of our construction, in Section~\ref{section:magnetic},
is a `magnetic'
modification of the symplectic
form on the standard filling $(DT^*\Sigma,\omegacan)$, by a method
reminiscent of the one used in \cite{mcdu91,geig95} to construct
symplectic manifolds with disconnected contact type boundaries.
This modification changes the symplectomorphism type of the
filling, but not its diffeomorphism type, and it makes the zero
section of the bundle $DT^*\Sigma$ symplectic. This
simple modification alone establishes the non-unique
fillability of $(ST^*\Sigma,\xican)$ up to symplectomorphism.

The zero section now being symplectic allows us to perform
a Gompf sum~\cite{gomp95,mcwo94} of $DT^*\Sigma$ with a copy of
$\Sigma\times\Sigma$
(equipped with the product symplectic form) along the
zero section and the diagonal, respectively; see
Section~\ref{section:Gompf-Sigma}.
In Section~\ref{section:pi1} we then use the rich topology of
$\Sigma\times\Sigma$ in order
to realise the desired fundamental group by further Gompf sums.
Finally, in Section~\ref{section:minimal} we show that our
fillings are minimal by construction.
\section{A `magnetic' filling}
\label{section:magnetic}
We continue to reserve the notation $\Sigma$ for a hyperbolic
surface, that is, a closed, oriented surface of genus $\geq 2$,
even though some parts of the following discussion apply
to all surfaces.

Write $\tlambda$ for the canonical $1$-form on the cotangent
bundle $T^*\Sigma$. In terms of the bundle projection
$\pi\co T^*\Sigma\rightarrow\Sigma$, this form is defined
as $\tlambda_{\bfu}=\bfu\circ T\pi$ for $\bfu\in T^*\Sigma$.
In local coordinates $q_1,q_2$ on $\Sigma$ and
dual coordinates $p_1,p_2$ on the fibres of $T^*\Sigma$ we have
$\tlambda=p_1\,\rmd q_1+p_2\,\rmd q_2=:\bfp\,\rmd\bfq$.

Given a $2$-form $\sigma$ on $\Sigma$, the corresponding
\emph{magnetic flow} is the Hamiltonian flow defined by the Hamiltonian
function $H=|\bfp|^2/2$ and the symplectic form $\rmd\tlambda+\pi^*\sigma$
on $T^*\Sigma$. Since the Euler class
of $T^*\Sigma$ is non-trivial, the Gysin sequence implies that
$\pi^*\sigma$ is exact when restricted to $T^*\Sigma\setminus\Sigma$.
This observation, which also applies to $T^*S^2$, was used in
\cite[Lemma~2.1]{beze15} to construct a symplectic form
on $DT^*S^2$ which near the boundary coincides with
$\rmd\tlambda=\omegacan$, and for which the zero section
$S^2\subset DT^*S^2$ is symplectic.

Here we present a simpler and more conceptual construction of such a
symplectic form on $DT^*\Sigma$. We equip $\Sigma$ with the
hyperbolic metric of constant curvature~$1$, and we write $J$ for
the associated complex structure, that is,
$J_{\bfq}$ defines the rotation of $T_{\bfq}\Sigma$
through an angle $\pi/2$. This gives rise to the pair of
Liouville--Cartan forms $\lambda,\mu$ on $ST^*\Sigma$.
The $1$-form $\lambda$ is simply the restriction of $\tlambda$ to
(the tangent spaces of) $ST^*\Sigma$.
The second $1$-form $\mu$ is defined by
$\mu_{\bfu}(\bft)=-\bfu(JT_{\bfu}\pi(\bft))$, where $\bfu\in ST^*\Sigma$
and $\bft\in T_{\bfu}(ST^*\Sigma)$. In local coordinates
$(\bfq,\bfp)$, this means that $\mu=p_1\,\rmd q_2-p_2\,\rmd q_1$.

With $\alpha$ denoting the connection $1$-form on $ST^*\Sigma$, we
then have
the structure equations
\begin{eqnarray*}
\rmd\lambda & = & \mu\wedge\alpha,\\
\rmd\mu     & = & \alpha\wedge\lambda,\\
\rmd\alpha  & = & -\lambda\wedge\mu;
\end{eqnarray*}
see \cite{sith67} or \cite{agz18}. The factor $-1$ in the third
equation corresponds to the choice of a constant curvature $-1$ metric.
We remark that $\lambda\wedge\mu=\pi^*\omega_{\Sigma}$ is the pull-back
of the area form $\omega_{\Sigma}$ on~$\Sigma$, so the curvature
form of $\alpha$ is $-\omega_{\Sigma}$, and $\alpha\wedge\lambda\wedge\mu$
is a volume form on $ST^*\Sigma$. The Euler class $e$ of the
bundle $ST^*\Sigma\rightarrow\Sigma$
is then given by $[\omega_{\Sigma}/2\pi]\in H^2(\Sigma;\Z)$, which by
Gau{\ss}--Bonnet means that $e=2g-2$, as it should.

The canonical $1$-form $\tlambda=\bfp\,\rmd\bfq$ is homogeneous of degree $1$
in the fibre coordinates. Thus, writing $r$ for the radial coordinate
on $T^*\Sigma$, and identifying $\lambda$ with its pull-back to
$T^*\Sigma\setminus\Sigma$ under radial projection, we have
$\tlambda=r\lambda$ on $T^*\Sigma\setminus\Sigma$. Likewise,
we regard $\alpha$ as a $1$-form on $T^*\Sigma\setminus\Sigma$.
Notice that $\alpha$ does not extend smoothly into the zero section,
though $r^2\alpha$ does, since fibrewise $\alpha$ is simply the
angular form. The $2$-form $\rmd\alpha=-\pi^*\omega_{\Sigma}$,
however, is a horizontal form defined on all of $T^*\Sigma$.

\begin{prop}
\label{prop:magnetic}
Let $f=f(r)$ be a smooth function $DT^*\Sigma\rightarrow [0,1]$ that is
identically $1$ near $r=0$ and identically zero near $r=1$. Then,
for $\varepsilon>0$ sufficiently small, the $2$-form
\[\omegamag:=\omegacan+\varepsilon\rmd(f\alpha)\]
is a (non-exact!)\ symplectic form on $DT^*\Sigma$ that coincides
with $\omegacan$ near the boundary, and for which the zero section
$\Sigma\subset DT^*\Sigma$ is symplectic. In particular,
$(DT^*\Sigma,\omegamag)$ is a strong symplectic filling
of $(ST^*\Sigma,\xican)$ that is not symplectomorphic to
the standard filling.
\end{prop}

\begin{proof}
Near $r=0$, where $f\equiv 1$, we have $\omega=\rmd\tlambda-\varepsilon
\pi^*\omega_{\Sigma}$, which is clearly symplectic, as can be seen
from the local description of $\rmd\tlambda$ as $\rmd\bfp\wedge\rmd\bfq$.

For $r>0$, we have
\begin{eqnarray*}
\omega & = & \rmd(r\lambda)+\varepsilon\rmd(f\alpha)\\
       & = & \rmd r\wedge\lambda+r\mu\wedge\alpha+\varepsilon
             f'(r)\,\rmd r\wedge\alpha-\varepsilon f(r)\,\lambda\wedge\mu,
\end{eqnarray*}
whence
\[ \omega^2=2\bigl(r-\varepsilon^2f(r)f'(r)\bigr)\,
\rmd r\wedge\alpha\wedge\lambda\wedge\mu,\]
which is nowhere zero for $\varepsilon>0$ sufficiently small,
since $f'(r)$ vanishes near the zero section.

In the exact filling $(DT^*\Sigma,\rmd\tlambda)$ we have $\int_{\Sigma}
\rmd\tlambda=0$ for every embedded copy of the closed surface
$\Sigma$ in $DT^*\Sigma$, whereas in the magnetic filling
$(DT^*\Sigma,\omegamag)$ the zero section $0_{\Sigma}\subset DT^*\Sigma$
is symplectic, so the two fillings are not symplectomorphic.
\end{proof}
\section{The Gompf sum with $\Sigma\times\Sigma$}
\label{section:Gompf-Sigma}
We now write $\sigma:=\omega_{\Sigma}$ for the area form
on $\Sigma$ corresponding to the hyperbolic metric, and
we denote the diagonal in $\Sigma\times\Sigma$ by $\Delta_{\Sigma}$.
The map $\bfq\mapsto(\bfq,\bfq)$ then defines a symplectic embedding
\[ (\Sigma,2\sigma)\stackrel{\cong}{\longrightarrow}\Delta_{\Sigma}
\subset (\Sigma\times\Sigma,\sigma\oplus\sigma).\]

\begin{lem}
\label{lem:normal}
The normal bundle of $\Delta_{\Sigma}\subset
(\Sigma\times\Sigma,\sigma\oplus\sigma)$ is orientation-reversingly
diffeomorphic to the normal bundle of the zero section
$0_{\Sigma}\subset (DT^*\Sigma,\omegamag)$.
\end{lem}

\begin{proof}
In the symplectic manifold $(\Sigma\times\Sigma,\sigma\ominus\sigma)$
(notice the sign!),
the diagonal $\Delta_{\Sigma}$ is a Lagrange submanifold, so by the
Weinstein neighbourhood theorem a neighbourhood
of the diagonal is symplectomorphic to a neighbourhood
of the zero section $0_{\Sigma}\subset(T^*\Sigma,\omegacan)$.
Hence, the normal bundle of $\Delta_{\Sigma}\subset
(\Sigma\times\Sigma,\sigma\ominus\sigma)$ --- with bundle
orientation determined by an orientation of $\Sigma$ and the ambient
orientation given by the symplectic form --- is isomorphic
to $T^*\Sigma$, which has Euler class $-\chi(\Sigma)$. In
$(\Sigma\times\Sigma,\sigma\oplus\sigma)$, the ambient orientation
is the opposite one, so the normal bundle of the diagonal has Euler class
$\chi(\Sigma)$.

The magnetic symplectic form $\omegamag$ is obtained as a deformation
of $\omegacan$ through symplectic forms, so they define the same ambient
orientation. This means that the Euler class of the normal
bundle of $0_{\Sigma}\subset(DT^*\Sigma,\omegamag)$ equals $-\chi(\Sigma)$.
\end{proof}

After rescaling $\sigma\oplus\sigma$ by a positive constant, we may
assume that the symplectic submanifolds $\Delta_{\Sigma}\subset
(\Sigma\times\Sigma,\sigma\oplus\sigma)$ and
$0_{\Sigma}\subset (DT^*\Sigma,\omegamag)$ are symplectomorphic.
With the preceding lemma we are then in the position to form the Gompf sum
\[ W_0:=(\Sigma\times\Sigma)\#_{\Delta_{\Sigma}=0_{\Sigma}} DT^*\Sigma.\]

Since $DT^*\Sigma$ with a (disc bundle) neighbourhood
$U$ of $0_{\Sigma}$ removed is simply
a collar on~$\partial U$, the symplectic $4$-manifold $W_0$
is diffeomorphic to $(\Sigma\times\Sigma)\setminus\Delta_{\Sigma}$.
Hence $W_0$ admits a locally trivial fibration over $\Sigma$
with fibre $\Sigma\setminus\{*\}$, that is, the base is aspherical,
and the fibre has the homotopy type of a bouquet of circles.
From the homotopy exact sequence of the fibration we infer that
$\pi_2(W_0)=0$. Thus, we have constructed a new symplectically
aspherical strong filling of $(ST^*\Sigma,\xican)$.
\section{Realising arbitrary fundamental groups}
\label{section:pi1}
Our aim now is to show that by performing further Gompf sums
on the symplectic $4$-manifold just constructed, we can realise any
finitely presented group as the fundamental group, as claimed
in Theorem~\ref{thm:main}.
\subsection{The idea of the construction}
Let $G$ be any finitely presentable group, and choose a presentation
\[ G=\langle g_1,\ldots,g_k\,|\,r_1,\ldots, r_{\ell}\rangle.\]
In \cite[Theorem~4.1]{gomp95}, Gompf constructs a symplectic $4$-manifold
$W_G$ with fundamental group $\pi_1(W_G)\cong G$. We shall construct a
strong symplectic filling of $(ST^*\Sigma,\xican)$ by first forming
the Gompf sum
\[ W_1:=W_G\#_T(\Sigma\times\Sigma)\#_{\Delta_{\Sigma}=0_{\Sigma}}
DT^*\Sigma\]
of $W_G$ with the filling $W_0$ found in the preceding section along
a suitable symplectic $2$-torus $T$ in the two constituents, and then
killing superfluous elements in the fundamental group by further
Gompf sums along symplectic tori with the elliptic surface
\[ V:=\CP^2\#9\overline{\CP^2}.\]
In order to understand the effect of these Gompf sums on the
fundamental group, and to fix the
notation, we need to recall a little more explicitly how $W_G$ is built.
\subsection{The Gompf manifold $W_G$}
Gompf starts with the symplectic $4$-manifold $\Sigma_k\times T^2$,
equipped with a product symplectic form, where $k$ is the number of
generators in the presentation of~$G$. Denote the standard generators
of $\pi_1(\Sigma_k)$ by $\alpha_1,\beta_1,\ldots,\alpha_k,\beta_k$.
We use the same notation for simple closed smooth curves in $\Sigma_k$
representing these generators. We choose these curves such that
$\alpha_i$ intersects $\beta_i$ in a single point, and there are
no further points of intersection. Similarly, we choose generators
$\alpha,\beta$ of $\pi_1(T^2)$.

The quotient group
\[ \pi_1(\Sigma_k)/\langle\beta_1,\ldots,\beta_k\rangle\cong
\langle\alpha_1,\ldots,\alpha_k\rangle\]
is the free group in $k$ generators. Choose immersed closed
curves $\gamma_1,\ldots,\gamma_{\ell}$ in $\Sigma_k$ representing
the relations $r_1,\ldots,r_{\ell}$ (under the identification
of $g_i$ with $\alpha_i$) in this free group. For ease of notation we
set $\gamma_{j+\ell}:=\beta_j$, $j=1,\ldots,k$. Then
\[ G\cong\pi_1(\Sigma_k)/\langle\gamma_1,\ldots,\gamma_{k+\ell}\rangle.\]

The $T_i:=\gamma_i\times\alpha$, $i=1,\ldots,k+\ell$, are immersed
Lagrange tori in $\Sigma_k\times T^2$. Further, with $z$ a point
in $\Sigma_k$ disjoint from the $\alpha$- and $\gamma$-curves,
we consider the symplectic torus $\{z\}\times T^2$.

Possibly after increasing the genus of $\Sigma_k$ (and correspondingly
extending the system of curves $\gamma_i$ such that we still have a
quotient description of $G$ as above), and after a small perturbation of the
product symplectic form as well as the~$T_i$, one  may assume them to be
disjointly embedded \emph{symplectic} tori with trivial normal bundle.
The separation of the tori is achieved by choosing the
$\gamma_i$ in disjoint $C^1$-close copies of $\Sigma$ inside
the $3$-manifold $\Sigma\times\beta$.

Likewise, the elliptic surface $V$ contains a symplectic torus
(a generic fibre) with trivial normal bundle and simply connected
complement. Then $W_G$ is obtained by performing the Gompf sum of
$\Sigma_k\times T^2$ with $k+\ell+1$ copies of~$V$, along
the $T_i$ and $\{z\}\times T^2$ in $\Sigma_k\times T^2$, and along
a generic fibre in each copy of~$V$. This has the effect of killing
the classes of $\gamma_1,\ldots,\gamma_{k+\ell},\alpha,\beta$ in
$\pi_1(\Sigma_k\times T^2)$.
\subsection{The Gompf sum of $W_G$ with the filling}
We now provide details of the Gompf sum
\[ W_1=W_G\#_T(\Sigma\times\Sigma)\#_{\Delta_{\Sigma}=0_{\Sigma}} DT^*\Sigma.\]
Recall that in the construction of $W_G$ one had to perform a Gompf
sum of $\Sigma_k\times T^2$ with the elliptic surface $V$ along
a torus $\{z\}\times T^2$ and a generic fibre, respectively.
We now choose a torus $\{z'\}\times T^2\subset\Sigma_k\times T^2$,
with $z'$ in the connected component of $z$ in
\[ \Sigma_k\setminus
\{\alpha_1,\ldots,\alpha_k,\gamma_1,\ldots,\gamma_{k+\ell}\}.\]

We denote the standard generators of $\pi_1(\Sigma)$ and
corresponding curves in $\Sigma$, where $\Sigma=\Sigma_g$,
$g\geq 2$, by $\lambda_1,\mu_1,\ldots,\lambda_g,\mu_g$, with the same
convention as above for the $\alpha_i,\beta_i$. In a second copy
of $\Sigma$, which we use to form the product $\Sigma\times\Sigma$,
we choose parallel curves $\lambda_1',\mu_1',\ldots,\lambda_g',\mu_g'$,
that is, with $\lambda_i'$ and $\mu_i'$ intersecting $\mu_i$
and $\lambda_i$, respectively, but none of the other curves.
This ensures that the tori $\lambda_i\times\lambda_i'$
and $\mu_i\times\mu_{i+1}'$ in $\Sigma\times\Sigma$
(where we read $\mu_{g+1}'$ as $\mu_1'$)
are pairwise disjoint, and disjoint from the diagonal $\Delta_{\Sigma}$.

These $2g$ tori have trivial normal bundle. This is clear topologically,
and also a consequence of the
fact that they are Lagrangian with respect to the symplectic form
$\sigma\oplus\sigma$ on $\Sigma\times\Sigma$; they are still Lagrangian
in $W_0=(\Sigma\times\Sigma)\#_{\Delta_{\Sigma}=0_{\Sigma}}DT^*\Sigma$
if the Gompf sum is performed in a sufficiently small neighbourhood of
$\Delta_{\Sigma}\subset\Sigma\times\Sigma$.
As we shall explain presently, by a small perturbation of the
symplectic form on $W_1$ one can make these tori symplectic of equal area.
The fibre connected sum
$\#_T$ then stands for the Gompf sum along $\{z'\}\times T^2\subset W_G$
and $\lambda_1\times\lambda_1'\subset W_0$.
\subsection{Surgering out superfluous generators}
Finally, we want to perform further Gompf sums on 
$W_1=W_G\#_T(\Sigma\times\Sigma)\#_{\Delta_{\Sigma}=0_{\Sigma}} DT^*\Sigma$
with copies of the elliptic surface $V$ so as to kill all
generators of $\pi_1$ coming from
$(\Sigma\times\Sigma\setminus\Delta_{\Sigma})$.

\begin{lem}
After a small perturbation of the symplectic form on $W_1$,
we may assume that the (previously Lagrangian) tori
\[ \lambda_1\times\lambda_1',\ldots,
\lambda_g\times\lambda_g',\mu_1\times\mu_2',\ldots,\mu_g\times\mu_1'\]
are symplectic of equal area.
\end{lem}

\begin{proof}
According to \cite[Lemma~1.6]{gomp95}, such a perturbation of the
symplectic form is possible whenever the homology classes of the
Lagrangian tori in $H_2(W_1;\R)$ lie in an affine subspace
that does not contain~$0$. This affine subspace consists of the
linear combinations where the coefficients sum up to~$1$, so this necessary
condition is certainly satisfied if the homology classes are linearly
independent.

To see that they are, suppose that
\[ a_1[\lambda_1\times\lambda_1']+\cdots+
a_g[\lambda_g\times\lambda_g']+b_1[\mu_1\times\mu_2']+\cdots+
b_g[\mu_g\times\mu_1']=0\in H_2(W_1;\R).\]
Let $\ell_1,m_1,\ldots,\ell_g,m_g$ and $\ell_1',m_1',\ldots,\ell_g',m_g'$
be $1$-forms that constitute bases of $H^1(\Sigma;\R)$ dual
to $\lambda_1,\mu_1,\ldots,\lambda_g,\mu_g$ and
$\lambda_1',\mu_1',\ldots,\lambda_g',\mu_g'$, respectively.
We use the same notation for the pull-backs of these $1$-forms first
to $\Sigma\times\Sigma$ under the projection to the first or second
factor, respectively, and then to
$(\Sigma\times\Sigma)\setminus\Delta_{\Sigma}$ under the inclusion.

The evaluation of one of the $2$-forms
\[ \ell_1\wedge\ell_1',\ldots,\ell_g\wedge\ell_g',
m_1\wedge m_2',\ldots,m_g\wedge m_1'\]
on the linear combination of homology classes above picks out exactly one
of the coefficients, so all coefficients have to be zero.
\end{proof}

\begin{lem}
The symplectic $4$-manifold $W_2$ obtained from $W_1$
by Gompf summing $2g-1$ copies of $V$ along the symplectic tori
from the preceding lemma --- except $\lambda_1\times\lambda_1'$, which was used
to sum $W_G$ and $W_0$ --- has fundamental group $\pi_1(W_2)\cong G$.
\end{lem}

\begin{proof}
Since the complement of a generic fibre in the elliptic surface $V$ is
simply connected, these $2g-1$ Gompf sums have the effect of killing
all generators coming from $\Sigma\times\Sigma$, except $\lambda_1$
and~$\lambda_1'$. It remains to show that the Gompf sum
$\#_T$ used in the preceding section for the construction of
$W_1$ kills $\lambda_1$ and $\lambda_1'$, while
leaving the contribution of $W_G$ to the fundamental group
unaffected.

The torus $\lambda_1\times\lambda_1'\subset(\Sigma\times\Sigma)
\setminus\Delta_{\Sigma}$ is homotopic in $W_1$
to a parallel copy of $\{z\}\times T^2\subset\Sigma_k\times T^2$,
and hence further to a fibre of the elliptic surface $V$
(in the complement of the fibre used for the Gompf sum with
$\Sigma_k\times T^2$ along $\{z\}\times T^2$). It follows that
$\lambda_1$ and $\lambda_1'$ are trivial in $\pi_1(W_1)$.

However, a meridional loop around $\{z'\}\times T^2\equiv
\lambda_1\times\lambda_1'$ created by the Gompf sum $\#_T$ might
be homotopically non-trivial. As a representative of this
homotopy element we may take a small loop $\gamma'$
in $\Sigma_k$ around $z'\in\Sigma_k$.

Let $\gamma$ be a small loop around $z\in\Sigma_k$, oriented in the
same way as~$\gamma'$. We may assume that the curves
$\gamma_1,\ldots,\gamma_{k+\ell}$ used in the construction of $W_G$
lie in other copies of $\Sigma_k$ in $\Sigma_k\times\beta$ than
$\gamma$ and~$\gamma'$. Then $\gamma'$ is homotopic in
$\Sigma_k\setminus\{z,z'\}$ (and hence in $W_1$) to the concatenation
$\prod_j[\alpha_j,\beta_j]\gamma^{-1}$.

Now, $\gamma$ becomes homotopically trivial in the elliptic surface
$V$ summed along $\{z\}\times T^2$ when we build~$W_G$. Likewise
the $\beta_j=\gamma_{j+\ell}$ become homotopically trivial
thanks to the Gompf sums with $V$ along $\beta_j\times\alpha$. It follows that
\[ \gamma'\simeq\prod_j[\alpha_j,\beta_j]\gamma^{-1}\simeq *\]
after these Gompf sums, so $\gamma'$ does not contribute non-trivially
to the fundamental group of~$W_1$.
\end{proof}

This concludes the proof of Theorem~\ref{thm:main}.
\section{Minimality of the filling}
\label{section:minimal}
Many classification results in symplectic topology hinge on
having a good moduli theory of $J$-holomorphic curves.
In particular, one typically wants to work in a setting
where the compactness of these moduli spaces can be
guaranteed. In higher dimensions, one often needs to require
symplectic asphericity for this to be the case;
in dimension~$4$, where one can appeal to positivity
of intersections, minimality of the symplectic manifolds in question
usually suffices. The results about uniqueness of
fillings of the standard contact spheres cited in the
introduction are an instance of these phenomena; see
also \cite[Section 4.5]{mcsa17}.

In this section we show that the fillings $W_2$ found
in the proof of Theorem~\ref{thm:main} are minimal by
construction, so there is no need to blow down
symplectic $(-1)$-spheres to make the fillings minimal.

We briefly recall how the elliptic fibration
on $V=\CP^2\# 9\overline{\CP^2}$ is constructed.
One starts with a generic pencil of cubic curves
in $\CP^2$. After blowing up the nine points where
the curves intersect, they form the fibres of an
elliptic fibration $V\rightarrow\CP^1$; the projection
to $\CP^1$ is given by mapping each elliptic curve
to the point of intersection with one of the nine exceptional
spheres created by the blow-ups.

\begin{lem}
\label{lem:VFminimal}
The complement $V\setminus F$ of a generic (non-singular)
fibre $F\subset V$ is minimal, but not symplectically
aspherical.
\end{lem}

\begin{proof}
First we show that $V\setminus F$ is not symplectically
spherical. The manifold $V\setminus F$ is simply connected;
a null-homotopy of the meridional loop around $F$
is provided by the transverse intersection of $F$
with one of the exceptional spheres. It follows that the
Hurewicz homomorphism $\pi_2(V\setminus F)\rightarrow
H_2(V\setminus F)$ is an isomorphism.
So the class $[F']\in H_2(V\setminus F)$ represented
by a generic fibre $F'\neq F$ is spherical, and any
such fibre is a symplectic submanifold.

Regarding minimality, we show something a little stronger:
there is no embedded surface $S\subset V\setminus F$
(symplectic or not) of self-intersection~$-1$. 

Write $A\in H_2(V)$ for the class of self-intersection $1$ represented
by $\CP^1\subset\CP^2$ in the complement of the nine
points that are blown up. The exceptional spheres
created by the blow-up represent classes $E_1,\ldots,E_9
\in H_2(V)$ of self-intersection~$-1$.

Now let $S\subset V$ be a surface embedded in the complement
of~$F$. Identifying $S$ with the homology class it
represents, we write
\[ S=aA+\sum_{j=1}^9e_jE_j,\]
with $a,e_1,\ldots,e_9\in\Z$. With the same convention for~$F$,
we have the intersection numbers
\[ F\bullet A=3\;\;\;\text{and}\;\;\; F\bullet E_j=1,\]
since $F$ comes from a cubic curve in $\CP^2$, the fibre $F$
intersects each exceptional sphere transversely in a single
point, and all these surfaces are symplectic.

The requirement that $S$ be disjoint from $F$ then translates into
\[ 0=F\bullet S=3a+\sum_{j=1}^9e_j.\]
Then the self-intersection number
\begin{eqnarray*}
S\bullet S & = & a^2-\sum_j e_j^2\\
           & = & (a+\sum_j e_j)\cdot(a-\sum_j e_j)
                 +2\sum_{i<k} e_ie_k\\
           & = & -2a\cdot 4a+2\sum_{i<k}e_ie_k
\end{eqnarray*}
is even.
\end{proof}

It follows that whenever we used a Gompf sum with $V$
along $F$ in the construction of the filling, which will
be the case for `most' fundamental groups, the filling
will not be symplectically aspherical. Minimality of the
filling, however, follows from \cite[Theorem~1.1]{ushe06}.
This theorem of Usher's states that with a few obvious exceptions,
the Gompf sum of two closed symplectic $4$-manifolds is minimal.
The exceptions are: (i) when there is a symplectic
$(-1)$-sphere in one of the two summands, disjoint from
the surface used for the fibre sum; (ii) if one summand
admits the structure of an $S^2$-bundle over a surface,
a section of the bundle is used for the Gompf sum, and the
second summand is not minimal.

Usher's theorem relies on Gromov--Witten theory and applies
to closed symplectic manifolds only (unless one wishes to appeal
to Seiberg--Witten theory), so we cannot use it directly
to show that $W_2$ is minimal. However, the magnetic deformation
in the proof of Proposition~\ref{prop:magnetic} can be localised
near the zero section of $DT^*\Sigma$. So by the Weinstein
neighbourhood theorem we may assume that the deformation of the symplectic
form is compactly supported in a neighbourhood symplectomorphic to
a neighbourhood of the diagonal $\Delta_{\Sigma}$ in $(\Sigma\times\Sigma,
\sigma\ominus\sigma)$; cf.\ the proof of Lemma~\ref{lem:normal}.

Then the Gompf sum of two copies of $\Sigma\times\Sigma$ along
the diagonal~$\Delta_{\Sigma}$, one copy equipped with the symplectic form
$\sigma\oplus\sigma$, the other, with the magnetic deformation
of $\sigma\ominus\sigma$, defines a closed symplectic manifold
$\widehat{W}_0$ into which $W_0$ embeds symplectically. Likewise,
$W_i$ embeds into a closed symplectic manifold $\widehat{W}_i$, $i=1,2$.

Usher's theorem shows that $\widehat{W}_2$, and hence $W_2$, is minimal.
Case (ii) of Usher's theorem does not occur
in our construction; case (i) is ruled out
by Lemma~\ref{lem:VFminimal} and the fact that the
other summands are topologically aspherical.


\begin{thebibliography}{10}
%
\bibitem{agz18}
\textsc{P. Albers, H. Geiges and K. Zehmisch},
Reeb dynamics inspired by Katok's example in Finsler geometry,
\textit{Math. Ann.}
\textbf{370} (2018), 1883--1907.
%
\bibitem{bgz19}
\textsc{K. Barth, H. Geiges and K. Zehmisch},
The diffeomorphism type of symplectic fillings,
\textit{J. Symplectic Geom.}
\textbf{17} (2019), 929--971.
%
\bibitem{beze15}
\textsc{G. Benedetti and K. Zehmisch},
On the existence of periodic orbits for magnetic systems on the two-sphere,
\textit{J. Mod. Dyn.}
\textbf{9} (2015), 141--146.
%
\bibitem{elia04}
\textsc{Ya. Eliashberg},
A few remarks about symplectic filling,
\textit{Geom. Topol.}
\textbf{8} (2004), 277--293.
%
\bibitem{etny04}
\textsc{J. B. Etnyre},
On symplectic fillings,
\textit{Algebr. Geom. Topol.}
\textbf{4} (2004), 73--80.
%
\bibitem{etho02}
\textsc{J. B. Etnyre and K. Honda},
On symplectic cobordisms,
\textit{Math. Ann.}
\textbf{323} (2002), 31--39.
%
\bibitem{geig95}
\textsc{H. Geiges},
Examples of symplectic $4$-manifolds with disconnected boundary of
contact type,
\textit{Bull. London Math. Soc.}
\textbf{27} (1995), 278--280.
%
\bibitem{geig06}
\textsc{H. Geiges},
Contact Dehn surgery, symplectic fillings, and Property P for knots,
\textit{Expo. Math.}
\textbf{24} (2006), 273--280.
%
\bibitem{geig08}
\textsc{H. Geiges},
\textit{An Introduction to Contact Topology},
Cambridge Studies in Advanced Mathematics \textbf{109},
Cambridge University Press, Cambridge (2008).
%
\bibitem{gkz23}
\textsc{H. Geiges, M. Kwon and K. Zehmisch},
Diffeomorphism type of symplectic fillings of unit
cotangent bundles,
\textit{J. Topol. Anal.}
\textbf{15} (2023), 683--705.
%
\bibitem{gomp95}
\textsc{R. E. Gompf},
A new construction of symplectic manifolds,
\textit{Ann. of Math. (2)}
%
% NOTE TO PRINTER: (2) refers to 2nd series and is part
% of the journal name. It should not be removed.
%
\textbf{142} (1995), 527--595.
%
\bibitem{grom85}
\textsc{M. Gromov},
Pseudoholomorphic curves in symplectic manifolds,
\textit{Invent. Math.}
\textbf{82} (1985), 307--347.
%
\bibitem{hind00}
\textsc{R. Hind},
Holomorphic filling of $\RP^3$,
\textit{Commun. Contemp. Math.}
\textbf{2} (2000), 349--363.
%
\bibitem{lmy17}
\textsc{T.-J. Li, C. Y. Mak and K. Yasui},
Calabi--Yau caps, uniruled caps and symplectic fillings,
\textit{Proc. Lond. Math. Soc. (3)}
%
% NOTE TO PRINTER: (3) refers to 3rd series and is part
% of the journal name. It should not be removed.
%
\textbf{114} (2017), 159--187.
%
\bibitem{mcwo94}
\textsc{J. D. McCarthy and J. G. Wolfson},
Symplectic normal connected sum,
\textit{Topology}
\textbf{33} (1994), 729--764.
%
\bibitem{mcdu90}
\textsc{D. McDuff},
The structure of rational and ruled symplectic $4$-manifolds,
\textit{J. Amer. Math. Soc.}
\textbf{3} (1990), 679--712.
%
\bibitem{mcdu91}
\textsc{D. McDuff},
Symplectic manifolds with contact type boundaries,
\textit{Invent. Math.}
\textbf{103} (1991), 651--671.
%
\bibitem{mcsa17}
\textsc{D. McDuff and D. Salamon},
\textit{Introduction to Symplectic Topology}
(3rd edn),
Oxford Graduate Texts in Mathematics \textbf{27}, Oxford University Press,
Oxford (2017).
%
%
\bibitem{sith67}
\textsc{I. M. Singer and J. A. Thorpe},
\textit{Lecture Notes on Elementary Topology and Geometry}
(Scott, Foresman and Co., Glenview, Ill., 1967).
%
\bibitem{siva17}
\textsc{S. Sivek and J. van Horn-Morris},
Fillings of unit cotangent bundles,
\textit{Math. Ann.}
\textbf{368} (2017), 1063--1080.
%
\bibitem{stip02}
\textsc{A. I. Stipsicz},
Gauge theory and Stein fillings of certain $3$-manifolds,
\textit{Turkish J. Math.}
\textbf{26} (2002), 115--130.
%
\bibitem{ushe06}
\textsc{M. Usher},
Minimality and symplectic sums,
\textit{Int.\ Math.\ Res.\ Not.}
\textbf{2006}, Article ID 49857, 17 pp.
%
\bibitem{wend10}
\textsc{C. Wendl},
Strongly fillable contact manifolds and $J$-holomorphic foliations,
\textit{Duke Math. J.}
\textbf{151} (2010), 337--384.
%
\end{thebibliography}
\end{document}